\documentclass[12pt]{amsart}

\usepackage{cancel}
\usepackage{latexsym, amssymb, amsmath,esint}
\usepackage{soul}
\usepackage{amsfonts, graphicx}
\usepackage{graphicx,color}

\usepackage{amsmath, amssymb, amsfonts, mathrsfs, mathtools}
\usepackage{latexsym, esint, graphicx, xcolor, color}
\usepackage{wasysym, stmaryrd}
\usepackage{hyperref}
\usepackage[alphabetic]{amsrefs}

\newcommand{\be}{\begin{equation}}
\newcommand{\ee}{\end{equation}}
\newcommand{\beq}{\begin{eqnarray}}
\newcommand{\eeq}{\end{eqnarray}}

\usepackage{wasysym,stmaryrd}
\newtheorem{thm}{Theorem}[section]

\newtheorem{lma}[thm]{Lemma}

\newtheorem{prop}[thm]{Proposition}
\newtheorem{cor}[thm]{Corollary}

\theoremstyle{remark}
\newtheorem{rem}[thm]{Remark}

\numberwithin{equation}{section}

\newtheorem{claim}{Claim}[section]

\def\be{\begin{equation}}
\def\ee{\end{equation}}
\def\bee{\begin{equation*}}
\def\eee{\end{equation*}}

\def\K{K\"ahler }
\mathtoolsset{showonlyrefs}

\def\Ric{\text{\rm Ric}}
\def\Rm{\text{\rm Rm}}

\def\e{\varepsilon}

\def\a{{\alpha}}
\def\b{{\beta}}

\begin{document}

\title[Volume growth and gap Theorem]{Volume growth and gap Theorem for manifolds with non-negative complex sectional curvature}

\author{Pak-Yeung Chan}
\address[Pak-Yeung Chan]{Department of Mathematics, National Tsing Hua University, Hsin-Chu, Taiwan}
\email{pychan@math.nthu.edu.tw}

\author{Man-Chun Lee}
\address[Man-Chun Lee]{Department of Mathematics, The Chinese University of Hong Kong, Shatin, Hong Kong, China}
\email{mclee@math.cuhk.edu.hk}

\subjclass[2020]{53C20, 53C21.}
\keywords{volume growth, gap Theorem, non-negative curvature}

\date{\today}

\begin{abstract}
The combined work of Liu, Ni, and Ni–Tam established that a complete K\"ahler manifold with positive bisectional curvature has Euclidean volume growth if and only if its curvature decays quadratically on average. In this work, we use the approximation scheme developed by Cabezas-Rivas and Wilking to establish the Riemannian analogue under the assumption of nonnegative complex sectional curvature. As an application, we prove an optimal gap theorem for manifolds with nonnegative complex sectional curvature. These results generalize the authors’ earlier work. 
\end{abstract}

\maketitle

\section{Introduction}
\label{sec: introduction}
Let $(M^n,g)$ be a complete Riemannian manifold. In this work, we study the geometric implications of the scale-invariant quantity
\begin{equation}
    k(x,r):=r^2\fint_{B_g(x,r)} \mathrm{scal}(g)\,d\mathrm{vol}_g,
\end{equation}
which measures the average quadratic decay of scalar curvature.
We are particularly interested in its influence on the large-scale asymptotic geometry of $(M,g)$ under nonnegative curvature assumptions. From an analytic viewpoint, this corresponds to the borderline regime for Morrey-type estimates. This quantity already plays a prominent role in determining the complex structure of a complete \K manifold with nonnegative bisectional curvature $\mathrm{BK}\geq 0$. To the best of our knowledge, the relationship between the decay rate of $k(x,r)$ and the existence of holomorphic functions of polynomial growth was first systematically investigated by Ni, Shi, and Tam \cite{NiShiTam2001}. It also appears in Shi's program \cite{Shi1997} for resolving Yau's uniformization conjecture using Ricci flow; see also \cite{ChenZhuTang}. Of particular interest is the case in which this scale-invariant quantity is uniformly bounded. In this setting, the combined work of Liu, Ni, and Ni--Tam established the following equivalence.

\begin{thm}[\cites{Ni2004,NiTam2013,Liu2016,Liu2018}]\label{thm:Kahler}
Let $(M^{n},g)$ be a complete noncompact \K manifold with $\mathrm{dim}_{\mathbb{C}}(M)=n\geq 2$ and quasi-positive bisectional curvature $\mathrm{BK}$. Then the following are equivalent:
\begin{enumerate}
    \item[(i)] $\displaystyle \mathrm{AVR}(g):=\lim_{r\to+\infty} \frac{\mathrm{Vol}_g(x_0,r)}{\omega_{2n}r^{2n}}>0$;
    \item[(ii)] There exists a nontrivial holomorphic function $f$ of polynomial growth;
    \item[(iii)] There exists a constant $C>0$ such that, for all $x\in M$ and $r>0$,
    \begin{equation}\label{eqn:k}
        k(x,r):=r^2\fint_{B_g(x,r)} \mathrm{scal}\,d\mathrm{vol}_g \leq C.
    \end{equation}
\end{enumerate}
\end{thm}

Under the assumption $\mathrm{BK}\geq 0$, Ni \cite{Ni2012} proved that if the curvature decay condition (iii) is strengthened to $k(x_0,r)=o(1)$ as $r\to+\infty$ for some $x_0\in M$, then the complete \K manifold must be flat. See also \cite{NiNiu2020} for an extension to a weaker curvature condition. We also refer to the recent work of Liu \cite{Liu2024}, which established a sharp relationship among these conditions under the assumption of Euclidean volume growth.

While this equivalence in \K geometry is now well understood, a precise Riemannian counterpart of Theorem~\ref{thm:Kahler} remains elusive in general. In Ricci curvature geometry, a related problem is to bound $k(x,r)$ uniformly from above in terms of local lower Ricci curvature bounds and volume noncollapsing assumptions; see \cites{JiangLiu,JiangNaber} for recent developments. Under the stronger assumption of nonnegative sectional curvature, Petrunin \cite{Petrunin2009} (see also \cite{Li}) showed that every complete Riemannian manifold $(M^n,g)$ with $\mathrm{sec}\geq 0$ satisfies
\begin{equation}\label{pet}
    \sup_{r>0} r^{2-n}\int_{B_g(x_0,r)}\mathrm{scal}(g)\,d\mathrm{vol}_g\leq C_n
\end{equation}
for every $x_0\in M$, where $C_n$ depends only on the dimension. Thus, under $\mathrm{sec}\geq 0$, Euclidean volume growth, as in (i), implies average quadratic curvature decay, as in (iii). The main difficulty is therefore to deduce volume noncollapsing from an appropriate curvature decay assumption. Intuitively, this asks whether regularity precludes collapse. Related ideas were explored in dimension three by Reiris \cite{Reiris} and Lott \cite{Lott}. In previous joint work with Li \cite{ChanLeeLi2026}, we used a level-set approach based on the Cheeger--Gromoll exhaustion \cite{CheegerGromoll1972} to prove a Riemannian analogue of Theorem~\ref{thm:Kahler} in dimension three under $\mathrm{sec}\geq 0$; see also \cite{Hu2026}.

In this work, we consider complex sectional curvature $\mathrm{sec}^{\mathbb{C}}$ in arbitrary dimensions. To clarify our notation, an algebraic curvature tensor $\mathrm{Rm}$ is said to have nonnegative complex sectional curvature, denoted by $\mathrm{sec}^{\mathbb{C}}\geq 0$, if, for every complex two-dimensional subspace $\Sigma$, its complexification satisfies $\Rm(u,v,\bar u,\bar v)\geq 0$ for every unitary basis $\{u,v\}$ of $\Sigma$. This condition is also denoted by $\mathrm{WPIC}2$ in the literature; we refer to \cite{BrendleBOOK} for its important role in the differentiable sphere theorem. In dimension three, it is equivalent to nonnegative sectional curvature $\mathrm{sec}\geq 0$. Our main result extends \cite[Corollary 1.3]{ChanLeeLi2026} to all dimensions $n\geq 3$ and provides a Riemannian counterpart to the equivalence between (i) and (iii) in Theorem~\ref{thm:Kahler}, formulated in terms of the universal cover. More precisely, we show that a complete noncompact manifold $(M^n,g)$ with $\mathrm{sec}^{\mathbb{C}}\geq 0$ and $\sup_{r>0}k(x_0,r)<+\infty$ for some $x_0\in M$ has a universal cover $(\tilde M^n,\tilde g)$ that is topologically Euclidean and satisfies $\mathrm{AVR}(\tilde g)>0$.

\begin{thm}\label{thm:main}
Let $(M^n,g)$ be a complete Riemannian manifold satisfying the following conditions:
\begin{enumerate}
    \item[(i)] $\mathrm{sec}^{\mathbb{C}}\geq 0$;
    \item[(ii)] There exist $x_0\in M$ and $C_0>0$ such that, for all $r>0$,
    \begin{equation}
        k(x_0,r):=r^2\fint_{B_g(x_0,r)} \mathrm{scal}\,d\mathrm{vol}_g\leq C_0.
    \end{equation}
\end{enumerate}
Then its universal cover $(\tilde M^n,\tilde g)$ is topologically Euclidean and has Euclidean volume growth.
\end{thm}

\begin{rem}
Under the assumption of Theorem \ref{thm:main}, if sectional curvature of $(M,g)$ is positive at a point, then it follows from \cite[Theorem 1.2]{CabezasWilking2015} or \cite{Perelman1994}  that $M$ is topologically Euclidean and hence $(M^n,g)$ has Euclidean volume growth. In combining the curvature estimate \eqref{pet} by Petrunin \cite{Petrunin2009}, it shows that in this case the maximal volume growth is equivalent to the boundedness of $\sup_{r>0}k(x_0,r)$ at a point $x_0\in M$. In this way, Theorem \ref{thm:main} can be viewed as the Riemannian analog of Theorem \ref{thm:Kahler}.
\end{rem}

\begin{rem} As discussed in \cite{ChanLeeLi2026}, the warped product metric
\[
dr^2+\tanh^2(r)\cdot g_{\mathbb S^{n-1}}
\]
on $\mathbb{R}^n$, has nonnegative curvature operator, linear volume growth and $$\limsup_{r\to+\infty}k(x_0,r)=+\infty.$$ 
\end{rem}

As an application, we use the volume noncollapsing obtained in Theorem~\ref{thm:main} to study a gap phenomenon for manifolds with $\mathrm{sec}^{\mathbb{C}}\geq 0$ and vanishing asymptotic average curvature, that is, $k(x_0,r)=o(1)$ as $r\to+\infty$, where $k(x_0,r)$ is defined in \eqref{eqn:k}. Gap theorems of this type in Riemannian category were first studied by the authors \cite{ChanLee2025}. Our aim is to establish a Riemannian analogue of Ni's optimal gap theorem \cite{Ni2012} (see also \cite{NiNiu2020}), which asserts that a complete noncompact \K manifold $(M^n,g)$ with $\mathrm{BK}\geq 0$ is flat if $k(x_0,r)=o(1)$ as $r\to+\infty$ for some $x_0\in M$. In this direction, we prove the following gap theorem, extending \cite[Theorem 1.1]{ChanLee2025}.

\begin{cor}\label{cor:gap}
Let $(M^n,g)$ be a complete Riemannian manifold satisfying the following conditions:
\begin{enumerate}
    \item[(i)] $\mathrm{sec}^{\mathbb{C}}\geq 0$;
    \item[(ii)] There exists $x_0\in M$ such that
    \begin{equation}
        k(x_0,r):=r^2\fint_{B_g(x_0,r)} \mathrm{scal}\,d\mathrm{vol}_g=o(1)
    \end{equation}
    as $r\to+\infty$.
\end{enumerate}
Then $(M,g)$ is flat.
\end{cor}

Our approach differs from that of \cite{ChanLeeLi2026} in its essential use of Ricci flow smoothing. This allows us to convert the asymptotic control of $k(x_0,r)$ into pointwise estimates and then apply the theorem of Gromoll and Meyer \cite{GromollMeyer1969} to obtain uniform volume noncollapsing. The main difficulty is to construct a Ricci flow with suitable quantitative estimates. Without a bounded curvature assumption, even short-time existence remains open in general. Under $\mathrm{sec}^{\mathbb{C}}\geq 0$, short-time existence of Ricci flow was first studied by Cabezas-Rivas and Wilking \cite{CabezasWilking2015}. Their construction relies on the Cheeger--Gromoll exhaustion \cites{CheegerGromoll1972,GromollMeyer1969} and the proof of the differentiable sphere theorem by Brendle and Schoen \cite{BrendleSchoen2009}. Although short-time existence was established, the resulting solution depends strongly on the local geometry of the initial manifold, and the construction generally provides no quantitative control. In this work, we refine and extend their construction under an asymptotic bound on $k(x_0,r)$. In particular, we show that the solution $g(t)$ constructed in \cite{CabezasWilking2015} exists for all time whenever the curvature satisfies an average bound centered at some point $x_0\in M$; see Theorem~\ref{thm:immort}.

\subsection*{Acknowledgments}  
P.-Y. Chan is supported by the Yushan Young Fellow Program of the Ministry of Education (MOE) Taiwan (MOE-108-YSFMS-0004-012-P1), and by NSTC grants 113-2115-M-007 -014 -MY2 and 115-2628-M-007 -006. M.-C. Lee is supported by Hong Kong RGC grants No. 14300623 and No. 14304225, and an Asian Young Scientist Fellowship.

\medskip

\subsection*{Disclosure on AI assistance.} The authors acknowledge the assistance of ChatGPT Pro, which suggested using a Busemann function to construct a space-time cut-off function. This construction played an important role in establishing local estimates for the approximating solutions.  The authors take full responsibility for the mathematical content of the paper.

\section{Reduction to universal cover with point soul}
In this section, we will discuss some basic properties about average curvature decay. We first show that the decay rate is independent of the center we considered.

\begin{lma}\label{lma:indep-centr} Suppose that $(M^n,g)$ is a complete manifolds with $\Ric\geq 0$.
For any non-negative smooth function $f:M\to \mathbb{R}$ and $x\in M$, if we denote 
\begin{equation}
k_f(x,r):=r^2\fint_{B_g(x,r)}f\,d\mathrm{vol}_g,
\end{equation}
then for any $x_1,x_0\in M$, we have 
\begin{equation}
k_f(x_1,r)\leq k_f(x_0,r+s)\cdot \left(\frac{r+2s}{r}\right)^n
\end{equation}
for all $r>0$, where $s:=d_g(x_0,x_1)>0$.
\end{lma}
\begin{proof}
Let $x_1\in M$ such that $d_{g}(x_1,x_0)=s>0$. Then for all $r>0$, 
\begin{equation}
B_g(x_1,r)\subseteq B_g(x_0,r+s)\subset B_g(x_1,r+2s)
\end{equation}
so that 
\begin{equation}
\begin{split}
k_f(x_1,r):&=r^2 \fint_{B_g(x_1,r)}f\,d\mathrm{vol}_g\\
&\leq \frac{r^2}{\mathrm{Vol}_g(x_1,r)} \int_{B_g(x_0,r+s)}f\,d\mathrm{vol}_g\\
&= k_f(x_0,r+s)\cdot \frac{r^2}{(r+s)^2}\cdot \frac{\mathrm{Vol}_g(x_0,r+s)}{\mathrm{Vol}_g(x_1,r)}\\
&\leq  k_f(x_0,r+s)\cdot \frac{r^2}{(r+s)^2}\cdot \frac{\mathrm{Vol}_g(x_1,r+2s)}{\mathrm{Vol}_g(x_1,r)}.
\end{split}
\end{equation}

Result follows by volume comparison.

\end{proof}

Next we show that the average curvature decay can be carried over to its universal cover.

\begin{lma}\label{lma:univer-k}
Let $(M^n,g)$ be a complete non-compact manifolds with $\Ric\geq 0$ and $\pi:(\tilde M^n,\tilde g)\to (M^n,g)$ be the covering projection from its universal cover. For any non-negative smooth function $f:M\to \mathbb{R}$, we have 
\begin{equation}
 \fint_{B_{\tilde g}(\tilde x_0,r)} f\circ \pi \,d\mathrm{vol}_{\tilde g} \leq 4^n  \fint_{B_g(x_0,r)}f \,d\mathrm{vol}_g
\end{equation}
for any $r>0$ and $x_0\in M$.
\end{lma}
\begin{proof}
For $x_0\in M$ and $r>0$, we fix a $\tilde x_0\in \pi^{-1}( x_0)$ and consider the finite set $A_r:=\{\gamma\in \Gamma: d_{\tilde g}(\tilde x_0,\gamma\tilde x_0)\leq r\}$ where $\Gamma $ denotes the deck group of $\pi$.

We let $s:M\to \tilde M$ be the lift. Suppose $\tilde y\in B_{\tilde g}(\tilde x_0,r) \setminus  \pi^{-1}\left(\mathrm{cut}_g(x_0) \right)$, then $y:=\pi(\tilde y)$ satisfies 
\begin{equation}
d_{g}(x_0,y)\leq d_{\tilde g}(\tilde x_0,\tilde y)<r.
\end{equation}

Since $\tilde y=\gamma (s(y))$ for some $\gamma\in \Gamma$,
\begin{equation}
\begin{split}
d_{\tilde g}(\tilde x_0,\gamma\tilde x_0)
&\leq d_{\tilde g}(\gamma\tilde x_0,\gamma(s(y)))+ d_{\tilde g}(\gamma(s(y)), \tilde x_0)\\
&\leq d_{\tilde g}(\tilde x_0,s(y))+ d_{\tilde g}(\tilde y, \tilde x_0)<2r.
\end{split}
\end{equation}

Therefore $B_{\tilde g}(\tilde x_0,r) \setminus  \pi^{-1}\left(\mathrm{cut}_g(x_0) \right)\subseteq \bigcup_{\gamma\in A_{2r}} s\left(B_g(x_0,r) \right)$, and hence 
\begin{equation}\label{eqn:cover-bdd-inte}
\begin{split}
\int_{B_{\tilde g}(\tilde x_0,r)} f\circ \pi\,d\mathrm{vol}_{\tilde g}&\leq \sum_{\gamma\in A_{2r}} \int_{\gamma D_r} f\circ \pi\,d\mathrm{vol}_{\tilde g}\\
&= |A_{2r}|\cdot  \int_{ D_r} f\circ \pi\,d\mathrm{vol}_{\tilde g}\\
&=|A_{2r}|\cdot  \int_{B_{g}(x_0,r)} f\,d\mathrm{vol}_{ g}\\
\end{split}
\end{equation}
where $D_r:=s\left(B_g(x_0,r) \right)$. 

It remains to estimate $N(s):=|A_{s}|$. By considering the deck transformation and fundamental domain of the covering, we have 
\begin{equation}\label{eqn:cover-s}
\begin{split}
\mathrm{Vol}_g(x_0,s)\cdot N(s)\leq \mathrm{Vol}_{\tilde g}(\tilde x_0,2s).
\end{split}
\end{equation}

By combining \eqref{eqn:cover-s} with \eqref{eqn:cover-bdd-inte}, we conclude using volume comparison that 
\begin{equation}
\begin{split}
\int_{B_{\tilde g}(\tilde x_0,r)} f\circ \pi\,d\mathrm{vol}_{\tilde g}&\leq N(2r)\cdot \int_{B_{g}(x_0,r)} f\,d\mathrm{vol}_{ g}\\
&\leq \frac{ \mathrm{Vol}_{\tilde g}(\tilde x_0,4r)}{\mathrm{Vol}_g(x_0,2r)}\cdot \int_{B_{g}(x_0,r)} f\,d\mathrm{vol}_{ g}\\
&\leq 4^n\frac{ \mathrm{Vol}_{\tilde g}(\tilde x_0,r)}{\mathrm{Vol}_g(x_0,r)}\cdot \int_{B_{g}(x_0,r)} f\,d\mathrm{vol}_{ g}.
\end{split}
\end{equation}
\end{proof}

Using Lemma~\ref{lma:univer-k}, we show that under some average curvature decay, the universal cover of a manifold with $\mathrm{sec}^\mathbb{C}\geq 0$ must be diffeomorphic to an Euclidean space. 
\begin{prop}\label{prop:soul-pt}
Suppose $(M^n,g)$ is a complete manifold such that 
\begin{enumerate}
\item[(i)] $\mathrm{sec}^\mathbb{C}\geq 0$;
\item[(ii)] there exists $x_0\in M$ such that 
$$k(x_0,r):=r^2\fint_{B_g(x_0,r)}\mathrm{scal}(g)\,d\mathrm{vol}_g =o(r^2)$$
as $r\to+\infty$.
\end{enumerate}
Then its universal cover $(\tilde M^n,\tilde g)$ of $(M^n,g)$ is diffeomorphic to $\mathbb{R}^n$.
\end{prop}
\begin{proof}

By \cite[Theorem 1.2]{CabezasWilking2015}, the universal cover $(\tilde M,\tilde g)$ of $(M,g)$ splits isometrically as $(\Sigma,g_\Sigma)\times (F,g_F)$ where $\Sigma^k$ is the $k$-dimensional soul of $M$ and $F^{n-k}$ is diffeomorphic to $\mathbb{R}^{n-k}$. 

We claim that $\Sigma$ must be a zero-dimensional so that the assertion follows. Suppose on the contrary, $\mathrm{dim}(\Sigma)>0$. We let $p\in \Sigma$ be the point such that $\pi\left((p,0^{n-k})\right)=x_0\in M$. We consider the ball $B_{\Sigma}(p,r)\Subset \Sigma$. By triangle inequality,
\begin{equation}
 B_{\tilde g}\left( (p,0^{n-k}) ,r\right)\subseteq B_{\Sigma}(p,r)\times B_F(0^{n-k},r)\subseteq B_{\tilde g}\left( (p,0^{n-k}) ,\sqrt{2}r\right)
\end{equation}
and hence the product metric formula implies $\mathrm{scal}(\tilde g)(x,y)=\mathrm{scal}(g_\Sigma)(x)+\mathrm{scal}(g_F)(y)\ge \mathrm{scal}(g_\Sigma)(x)$ for $(x,y)\in \Sigma\times F$,
\begin{equation}
\begin{split}
\int_{B_{\tilde g}\left( (p,0^{n-k}) ,\sqrt{2}r\right)}\mathrm{scal}(\tilde g)\,d\mathrm{vol}_{\tilde g}
&\geq  \int_{B_{\Sigma}(p,r)}\int_{B_F(0^{n-k},r)} \mathrm{scal}(\tilde g)\,d\mathrm{vol}_F \,d\mathrm{vol}_{\Sigma}\\
&\geq \mathrm{Vol}_{g_F}(0^{n-k},r)\cdot \int_{B_{\Sigma}(p,r)} \mathrm{scal}(g_\Sigma)\,d\mathrm{vol}_{\Sigma}\\
&\geq \mathrm{Vol}_{\tilde g}\left((p,0^{n-k}),r \right)\cdot \fint_{B_{\Sigma}(p,r)} \mathrm{scal}(g_\Sigma)\,d\mathrm{vol}_{\Sigma}.
\end{split}
\end{equation}

By Lemma~\ref{lma:univer-k}, $(ii)$ and the compactness of  $\Sigma$, the inequality for $r\geq \mathrm{diam}(\Sigma)$ shows that 
\begin{equation}
\fint_{\Sigma}\mathrm{scal}(g_\Sigma)\,d\mathrm{vol}_{\Sigma}=\fint_{B_{\Sigma}(p,r)} \mathrm{scal}(g_\Sigma)\,d\mathrm{vol}_{\Sigma}=o(1)
\end{equation}
as $r\to+\infty$. Since $\mathrm{sec}^\mathbb{C}(\Sigma)\geq 0$, this forces $g_\Sigma$ to be flat. But this contradicts with simply connectedness of $\tilde M$.
\end{proof}

\section{Compact approximation by Cabezas-Rivas and Wilking}  

Our approach to Theorem~\ref{thm:main} is based on the approximation scheme developed by Cabezas-Rivas and Wilking \cite{CabezasWilking2015}. By exploring the Cheeger–Gromoll exhaustion \cites{CheegerGromoll1972,GromollMeyer1969}, Cabezas-Rivas and Wilking constructed short-time complete solution of the Ricci flow starting from any complete manifold with $\mathrm{sec}^\mathbb{C}\geq 0$. We first recall the technical tool involved. We will use $g_0$ to denote the metric under consideration in Theorem~\ref{thm:main}.

\subsection{Approximation using Cheeger–Gromoll exhaustion }

Fix $x_0\in M$, we consider the set $\mathcal{R}$ consisting of all rays $\gamma:[0,+\infty)\to M$ with $\gamma(0)=x_0$. The busemann function $\hat b:M\to \mathbb{R}$ is defined to be 
\begin{equation}
\hat b(x):=\sup_{\gamma\in\mathcal{R}} \left\{ \lim_{s\to+\infty} \left(s-d_{g_0}(x,\gamma(s)) \right)\right\}.
\end{equation}

By the work of Cheeger-Gromoll-Meyer \cites{CheegerGromoll1972,GromollMeyer1969}, $\hat b$ is proper and  $\inf_M \hat b\leq 0$. For our convenience, we instead consider the translated Busemann function $\b:=\hat b+\tilde L$ for some large $\tilde L$. To ease the discussion, we will assume the soul to be a point. This will be sufficient for our purpose.

We denote the level set 
\begin{equation}
C_\ell := \left\{x\in M:   \b(x)\leq \ell \right\}\subseteq M.
\end{equation}

\begin{lma}\label{lma:prop-buse}
Suppose $(M^n,g_0)$ is complete non-compact with $\mathrm{sec}\geq 0$, then the following holds: 
\begin{enumerate}
\item[(i)] $\b$ is a convex function; 
\item[(ii)] Each $C_\ell$ is totally convex compact set;
\item[(iii)] $\mathrm{dim}(C_\ell)=n$ for $\ell>0$ and $M=\cup_{\ell>0} C_\ell$;
\item[(iv)] $C_\ell$ has the structure of an embedded sub-manifold of $M$ with smooth totally geodesic interior and (possibly non-smooth) boundary;
\item [(v)] $\b$ is $1$-Lipschitz on $M$ and satisfies $|\nabla \b|=1$ almost everywhere;
\item [(vi)] We might choose $\tilde L$ large enough such that
\begin{equation}
\frac12 d_{g_0}(\cdot,p_0) \leq \b(\cdot )\leq d_{g_0}(\cdot,p_0)+\tilde L.
\end{equation}
\end{enumerate}
\end{lma}
\begin{proof}
The properties (i)-(iv) are from \cites{CheegerGromoll1972,GromollMeyer1969}. The $1$-Lipschitz property (v) follows from definition while (vi) follows from Toponogov’s comparison theorem, see \cite[Lemma 1]{Drees1994}. 
\end{proof}

\medskip

The following is implicitly proved in \cite[Section 6]{CabezasWilking2015}, by considering the doubling $D(C_\ell)$ of $C_\ell$ where $\ell\to+\infty$. As in \cite{CabezasWilking2015}, the tubular neighborhood $\cup_{x\in A}\overline B_r(x)$ of any subset $A$ in $M$ is denoted by $B_r(A)$.
\begin{prop}\label{prop:CW-appr}
Suppose $(M^n,g_0)$ is a complete non-compact manifold with $\mathrm{sec}^\mathbb{C}\geq 0$ whose soul is a point $p_0$, then there exists a sequence of  compact manifold $(M^n_i,g_{i,0},p_i)$ which converges to $(M^n,g_0,p_0)$ in the smooth pointed Cheeger-Gromov sense as $i\to+\infty$ such that 
\begin{enumerate}
\item[(i)] each $M_i$ is diffeomorphic to $\mathbb{S}^n$;
\item[(ii)] $\mathrm{sec}^\mathbb{C}(g_i)>0$;
\item[(iii)] each $(M_i,g_{i,0})$ has $\mathbb{Z}_2$-symmetry;
\item[(iv)] the fixed set $N_i$ of $\mathbb{Z}_2$-symmetry in (iii) is a smooth separating hyper-surface;
\item[(v)] there is $\ell_i\to+\infty$ such that $b_i:=\ell_i-d_{g_{i,0}}(\cdot,N_i)$ is  non-negative, $1$-Lipschitz, convex outside $N_i$ and converges to the function $\b$ under the Cheeger-Gromov convergence.
\end{enumerate}
\end{prop}
\begin{proof}
This follows implicitly from \cite[Proposition 6.7]{CabezasWilking2015}. We first recall part of the argument there. 

We let $\e_i\to 0$ and $\ell_i\to+\infty$. We will omit the index $i$ in $\ell_i$ for convenience. Let $C_\ell$ be the totally convex compact set in Lemma~\ref{lma:prop-buse}.  By \cite[Proposition 6.3]{CabezasWilking2015}, there exists a sequence of smooth closed hypersurface $D_{\ell,k}$ given by doubling $D(G_k)$ of some graph $G_k$ in $B_1(C_\ell)\times [0,1]$ converging to $D(C_\ell)$ in the Gromov-Hausdorff sense as $k\to+\infty$. And \cite[Corollary 6.4 \& Proposition 6.6]{CabezasWilking2015} shows that each $D_{\ell,k}$ admits  a  Ricci flow $(D_{\ell,k},g_{\ell,k}(t)),t\in (0,T_\ell]$ such that as $k\to+\infty$, it converges in Cheeger-Gromov topology to a limiting Ricci flow $(D_{\ell,\infty},g_{\ell,\infty}(t)),t\in (0,T_\ell]$ where $D_{\ell,\infty}$ is homeomorphic to the double of $C_\ell$, which is homeomorphic to $\mathbb{S}^n$. It then follows from strong maximum principle \cite{BrendleSchoen2008} that $\mathrm{sec}^\mathbb{C}(g_{\ell,\infty}(t))>0$ for all $t\in (0,T_\ell]$, see \cite[Proposition 3.3]{CabezasWilking2015}.

Moreover as $t\to 0$, $(D_{\ell,\infty},g_{\ell,\infty}(t))$ converges in Gromov-Hausdorff sense to $D(C_\ell)$. We choose $t_i\in (0,T_\ell]$ sufficiently small such that the Gromov-Hausdorff distance between $D(C_\ell)$ and $(D_{\ell,\infty},g_{\ell,\infty}(t_i))$ is at most $\e_i>0$. We might assume $t_i\to 0^+$.  For this $t_i$, we choose $k_i$ sufficiently large so that $(D_{\ell,\infty},g_{\ell,\infty}(t_i))$ is $\e_i$-close to $(D_{\ell_i,k_i},g_{\ell_i,k_i}(t_i))$ in the Cheeger-Gromov topology. In particular, we might assume $\mathrm{sec}^\mathbb{C}(g_{\ell_i,k_i}(t_i))>0$ since the underlying manifold is compact. Furthermore for any $\delta>0$, we can view $C_{\ell-\delta}$ as subset of $D_{\ell,k}$ from the construction in \cite[Proposition 6.3]{CabezasWilking2015}. If we choose $(M_i,g_{i,0},p_i)=(D_{\ell_i,k_i},g_{\ell_i,k_i}(t_i),p_0)$, then the Shi's estimates \cite[Theorem C.4]{CabezasWilking2015} implies the smooth pointed Cheeger-Gromov convergence to $(M,g_0,p_0)$ as $i\to+\infty$. The properties (i)-(ii) follows directly. Since compact Ricci flow preserves the symmetric, $(M_i,g_{i,0})$ is also symmetric across $N_i:=\partial G_{k_i}$, and thus it proves (iii)-(iv).  

Finally since $N_i$ is the fixed smooth hypersurface of the $\mathbb{Z}_2$ symmetry, $N_i$ is totally geodesic with respect to $g_{i,0}$. The function $b_i:=\ell_i-d_{g_{i,0}}(\cdot, N_i)$ is therefore $1$-Lipschitz by triangle inequality. It is convex with respect to $g_{\ell,\infty}(t_\ell)$ outside $N_i$ by second variation formula since $N_i$ is totally geodesic, see the proof of \cite[Lemma 4.8]{CabezasWilking2015}. Since the Gromov-Hausdorff distance between $D(C_\ell)$ and $(D_{\ell_i,k_i},g_{\ell_i,k_i}(t_i))$ is at most $\e_i$, on the $\mathbb{Z}_2$ orbit containing $p_0$, 
\begin{equation}
    b_i(x)=\ell_i-d_{g_0}(\cdot, \partial C_{\ell_i})+o(1)
\end{equation}
as $i\to+\infty$, and hence proves the convergence of $b_i$ to $\b$, 
using \cite[Lemma A.4]{CabezasWilking2015}.
\end{proof}

\subsection{Ricci flowing of approximation}

We let $(M_i,g_{i,0})$ be the compact manifold from  Proposition~\ref{prop:CW-appr}. In \cite{CabezasWilking2015}, the Ricci flow $g_i(t)$ starting from $g_{i,0}$ is considered in order to construct global solution $g(t)$  starting from $g_0$.  The following records a characterization of its lifespan, see also \cite[Proposition 4.3]{CabezasWilking2015}.
\begin{lma}\label{lma:lifespan}
Suppose $g_i(t),t\in [0,T_i)$ is the maximal solution to the Ricci flow starting from $(M^n_i,g_{i,0})$, then 
\begin{equation}
\limsup_{t\to T_i^-} \mathrm{Vol}_{g_i(t)}(M_i)=0.
\end{equation}
\end{lma}
\begin{proof}
 By the work of Brendle-Schoen \cite{BrendleSchoen2009}, the Ricci flow $g_i(t)$ converges to a shrinking sphere as $t\to T_i^-$. Since it is a (spherical) point singularity, result follows. 
\end{proof}

In \cite{CabezasWilking2015}, $T_i$ was controlled uniformly from below by local geometry of $(M,g_0)$ nearby $p_0$ using Lemma~\ref{lma:lifespan} and an estimate of Petrunin \cite{Petrunin2009}. We follow the spirit and will estimate the lifespan of each $g_i(t)$ for  volume non-collapsing. We first modify the distance like function $b_i$ to function adapted to  Ricci flow $g_i(t)$.

\begin{lma}\label{lma:heat-flowed}
On each $(M_i,g_i(t)),t\in [0,T_i)$, the function
$$u_i(x,t):=\ell_i-d_{g_i(t)}(x,N_i)$$
satisfies 
\begin{enumerate}
    \item[(i)] $0\leq \partial_t u_i\leq \Delta_{g_i} u_i$;
    \item[(ii)] $|\nabla u_i|\leq 1$
\end{enumerate}
in the sense of distribution, on $M_i\setminus N_i$. Moreover $u_i(x,t)\geq b_i(x)$ for all $(x,t)\in M_i\times [0,T_i)$.
\end{lma}
\begin{proof}
The Lipschitz bound follows from definition. The monotonicity follows directly from $\Ric(g_i(t))\geq 0$. Finally the property (i) follows from first variation using Ricci flow and second variation formula using totally geodesic of $N_i$, see the proof of \cite[Lemma 4.8]{CabezasWilking2015}.
\end{proof}

\section{Long-time existence of Ricci flow constructed by Cabezas-Rivas and Wilking}\label{sec:long-time}

In this section, we show that the Ricci flow constructed by Cabezas-Rivas and Wilking \cite{CabezasWilking2015} can be made immortal if the average curvature stays uniformly bounded on any compact set.

This will be based on the local smoothing via considering the Ricci flow starting from $(M_i,g_{i,0})$, from Proposition~\ref{prop:CW-appr}.

\subsection{Estimate of lifespan in approximating Ricci flows}
We now show that under the assumption of Theorem~\ref{thm:soul}, the Ricci flow $g_i(t)$ starting from $(M_i,g_{i,0})$ constructed in Proposition~\ref{prop:CW-appr} is almost immortal for $i\to+\infty$.

We will mimic the estimate using heat kernel in \cite{ChanLee2025}. We fix a smooth non-increasing function $\varphi$ on $[0,+\infty)$ such that $\varphi=1$ on $[0,1]$, vanishes outside $[0,2]$ and satisfies $|\varphi'|^2+|\varphi''|\leq 10$. We let $\phi:=\varphi^m$ where $m$ is a large integer to be chosen.  For $R>0$ and $i\to+\infty$, we define
\begin{equation}
\Phi_{i,R,m}(x,t):=\phi(u_i(x,t)/R)
\end{equation}
on the connected component $\Omega_i$ of $M_i\setminus N_i$ that contains $p_i$, where $u_i$ is from Lemma~\ref{lma:heat-flowed}. Since $u_i\geq b_i$ where $b_i(x)\approx \ell_i\to+\infty$ near $N_i$, $\Phi_{i,R}$ vanishes nearby $N_i$. We then extend $\Phi_{i,R,m}$ identically zero to the other orbit. Hence the evolution inequalities of $u_i$ in Lemma~\ref{lma:heat-flowed} hold in the support of $\Phi_{i,R}$.

Using the choice of cut-off function $\Phi_{i,R,m}$, we  study the following integrals: 
\begin{equation}\label{eqn:defining-W-V}
\left\{
\begin{array}{ll}
\displaystyle V_{i,R,m}(t):=\int_{M_i} \Phi_{i,R,m}\,d\mathrm{vol}_{g_i(t)}\\[3mm]
\displaystyle W_{i,R,m}(t):=\int_{M_i} \Phi_{i,R,m}\cdot\mathrm{scal}(g_i)\,d\mathrm{vol}_{g_i(t)}. 
\end{array}
\right.
\end{equation}

The motivation of doing this is to examine the non-collapsing via $V_{i,R,m}$, while its derivative is controlled by $W_{i,R,m}$. Since $W_{i,R,m}$ is initially controlled by the initial curvature decay, one expect that the volume integral $V_{i,R,m}$ can be controlled in a quantitative  way.

We first compute the evolution inequality of the weighted volume. 
\begin{lma}\label{eqn:evo-V}
For any large $R>0$ there is $N\in \mathbb{N}$ such that for $i>N$,
\begin{equation}
\partial_t V_{i,R,m}\geq  -W_{i,R,m} -\frac{10^33^n m^2}{R^2} \cdot\mathrm{Vol}_{g_{i,0}}\left(p_i,R\right)
\end{equation}
on $(0,T_i)$.
\end{lma}
\begin{proof}
For $t\in (0,T_i)$, using $0\leq \partial_t u_i\leq \Delta_{g_i}u_i$ we have 
\begin{equation}
\begin{split}
\partial_t V_{i,R,m}&= -W_{i,R,m}+\int_{M_i}  \frac1R\phi' \partial_t u_i \,d\mathrm{vol}_{g_i}\\
&\geq  -W_{i,R,m}+\int_{M_i}  \frac1R\phi' \Delta_{g_i}u_i \,d\mathrm{vol}_{g_i}\\
&= -W_{i,R,m}-\int_{M_i}  \frac1{R^2}\phi'' |\nabla u_i|^2 \,d\mathrm{vol}_{g_i}.
\end{split}
\end{equation}

Recall that $\Omega_i$ is the connected component in $M_i\setminus N_i$ containing $p_i$. From Lemma~\ref{lma:heat-flowed}, $u_i(x,t)\geq b_i(x)$, this implies 
\begin{equation}
\{ x\in \Omega_i : u_i(x,t)\leq 2R\} \subseteq \{x\in \Omega_i: b_i(x)\leq 2R\}.
\end{equation}
In view of the convergence of $b_i$ to $\b$ from Proposition~\ref{prop:CW-appr}, we may assume that $b_i(p_i)\leq 2R$ for large $R$ and all large $i$. and $b_i>2R+1$ on $\partial B_{g_{i,0}}(p_i,3R)$. Next we show that 
\begin{equation}
    \{x\in \Omega_i: b_i(x)\leq 2R\}\subseteq B_{g_{i,0}}(p_i,3R)\subseteq B_{g_i(t)}(p_i,3R).
\end{equation}
Note that the second inclusion follows from the fact that $\Ric(g_i(t))\geq 0$. It remains justify the first one. Fix any $q\in \{x\in \Omega_i: b_i(x)\leq 2R\}$, since $N_i$ is totally geodesic and $g_{i,0}$ is $\mathbb Z_2$-symmetric, there exists a $g_{i,0}$-minimizing geodesic $\gamma$ joining $p_i$ and $q$ which doesn't intersect $N_i$. If $q\notin B_{g_{i,0}}(p_i,3R)$, then $\gamma$ intersects $\partial B_{g_{i,0}}(p_i,3R)$ at some point, say $q'$. By the convexity of $b_i$ from Proposition~\ref{prop:CW-appr}, $2R+1<b_i(q')\leq \max{\{b_i(p_i), b_i(q)\}}\leq 2R.$ Contradiction. This shows that inclusion required. To summarize, we have shown
\begin{equation}\label{eqn:set-inc}
\begin{split}
\{ x\in M_i: u_i(x,t)\leq 2R\} &\subseteq \{x\in M_i: b_i(x)\leq 2R\}\\
&\subseteq B_{g_{i,0}}(p_i,3R)\subseteq B_{g_i(t)}(p_i,3R).
\end{split}
\end{equation}

Together with $|\nabla u_i|\leq 1$ from Lemma~\ref{lma:heat-flowed} and the definition of $\phi$, we conclude that 
\begin{equation}
\begin{split}
\partial_t V_{i,R,m}&\geq -W_{i,R,m} -\frac{10^3 m^2}{R^2} \cdot\mathrm{Vol}_{g_{i,0}}\left(p_i,3R\right)
\end{split}
\end{equation}
where we have used $g_i(t)\leq g_{i,0}$ from $\Ric(g_i(t))\geq 0$.

\end{proof}

Similarly, we have  the evolution inequality for $W_{i,R,r}$. 
\begin{lma}\label{eqn:evo-W}
There exists $C_1(n)>0$ such that for any large $R>0$, there exists $N\in \mathbb{N}$ so that for $i>N$,
\begin{equation}
\partial_tW_{i,R,m}\leq C_1m^2 R^{-2+\frac{2(n-2)}{m}} \left(W_{i,R,m}\right)^{1-\frac2m}
\end{equation}
for $t\in (0,T_i)$.
\end{lma}
\begin{proof}
First we note that from  \cite[Lemma 3.2]{ChanLee2025} that 
\begin{equation}
\left(\partial_t -\Delta_{g_i}-\mathrm{scal}(g_i) \right)\mathrm{scal}(g_i)\leq 0
\end{equation}
so that 
\begin{equation}\label{eqn:Energy-1}
\begin{split}
\partial_t W_{i,R,m}&\leq \int_{M_i}  \left( \partial_t +\Delta_{g_i}\right)\Phi_{i,R,m} \cdot \mathrm{scal}(g_i)\,d\mathrm{vol}_{g_i}\\
&=\int_{M_i}\left( \frac{\phi'}{R}\left( \partial_t +\Delta_{g_i}\right)u_i +\frac{\phi''}{R^2} |\nabla u_i|^2\right)\cdot \mathrm{scal}(g_i)\,d\mathrm{vol}_{g_i}.
\end{split}
\end{equation}
Here we have used the fact that the adjoint operator of $\partial_t-\Delta_{g_i(t)}-\mathrm{scal}(g_i(t))$ is precisely given by $\partial_t+\Delta_{g_i(t)}$.

Since $\phi'\leq 0$ and $0\leq \partial_t u_i\leq \Delta_{g_i} u_i$ on $M_i\setminus N_i$ from Lemma~\ref{lma:heat-flowed}, the first term in the bracket is non-positive. Using also $|\nabla u_i|^2\leq 1$ and the H\"older inequality, \eqref{eqn:Energy-1} is reduced to 
\begin{equation}
\begin{split}
\partial_t W_{i,R,m}&\leq \frac{10^3m^2}{R^2}\int_{M_i} \phi^{1-\frac2m}\mathrm{scal}(g_i)\,d\mathrm{vol}_{g_i}\\
&\leq \frac{10^3m^2}{R^2} \left(W_{i,R,m}\right)^{1-\frac2m}\cdot \left(\int_{B_{g_i}(p_i,3R)}\mathrm{scal}(g_i)\,d\mathrm{vol}_{g_i}\right)^\frac2m
\end{split}
\end{equation}
where we have used \eqref{eqn:set-inc} so that $\Phi_{i,R,m}$ is compactly supported in $B_{g_i(t)}(p_i,3R)$. The desired result now follows from  the rescaled version of  Petrunin's estimate \cite{Petrunin2009}, see \cite[Lemma 7.1]{CabezasWilking2015}, and $\mathrm{sec}(g_i)\geq 0$.
\end{proof}

\medskip
Finally, we remark that the initial value of $V_{i,R,m}$ and $W_{i,R,m}$ can be estimated in term of volume of balls in $(M,g_0)$.
\begin{lma}\label{lma:vol-init}
For $R$ sufficiently large, we have 
\begin{equation}
\lim_{i\to+\infty} V_{i,R,m}(0)\geq 2^{-n}\mathrm{Vol}_{g_0}(p_0,R).
\end{equation}
\end{lma}
\begin{proof}

For any given $R$ large, thanks to the convergence from Proposition~\ref{prop:CW-appr}  we have 
\begin{equation}
\begin{split}
\lim_{i\to+\infty} V_{i,R,m}(0)&=\int_M \phi\left(\frac{ \b(x)}{R}\right)  d\mathrm{vol}_{g_0}\\
&\geq \int_{ \b\leq R} d\mathrm{vol}_{g_0}.
\end{split}
\end{equation}

By Lemma~\ref{lma:prop-buse}, $\{ x\in M: \b\leq R\}$ contains $B_{g_0}(p_0,R/2)$ if $R$ is large. Hence volume comparison yields
\begin{equation}
\begin{split}
\lim_{i\to+\infty} V_{i,R,m}(0)
&\geq 2^{-n}\mathrm{Vol}_{g_0}(p_0,R).
\end{split}
\end{equation}
\end{proof}

Similarly, we have an analogous lower bound for $W_{i,R,m}$. 
\begin{lma}\label{lma:init-W}
For $R$ sufficiently large, we have
\begin{equation}
\lim_{i\to +\infty} W_{i,R,m}(0)\leq 4^{n-2}R^{-2} k(p_0,4R)\cdot \mathrm{Vol}_{g_0}(p_0,R).
\end{equation}
\end{lma}
\begin{proof}

For any given large $R$, thanks to the convergence from Proposition~\ref{prop:CW-appr} and Lemma~\ref{lma:prop-buse}, as well as volume comparison, we have 
\begin{equation}
\begin{split}
\lim_{i\to +\infty} W_{i,R}(0)&=\int_{\{\b(x)\leq 2R\}} \phi\left(\frac{\b(x)}{R} \right)\cdot\mathrm{scal}_{g_0}\,d\mathrm{vol}_{g_0}\\
&\leq \int_{B_{g_0}(p_0,4R)}   \mathrm{scal}_{g_0}\,d\mathrm{vol}_{g_0}\\
&= 4^{-2}R^{-2} k(p_0,4R)\cdot \mathrm{Vol}_{g_0}(p_0,4R)\\
&\leq 4^{n-2}R^{-2} k(p_0,4R)\cdot \mathrm{Vol}_{g_0}(p_0,R).
\end{split}
\end{equation}
\end{proof}

\medskip

Now we claim that the Ricci flow $(M_i,g_i(t)),t\in [0,T_i)$ is almost immortal as $i\to+\infty$. 
\begin{prop}\label{prop:almost-immort}
Suppose $(M^n,g_0)$ is a complete non-compact manifold with $\mathrm{sec}^\mathbb{C}\geq 0$, soul point $p_0$ and 
\begin{equation}
 \lim_{r\to+\infty}r^{-2}k(p_0,r)=0,\;\;\text{as}\;r\to+\infty.
\end{equation}
 If $g_i(t),t\in [0,T_i)$ is the maximal solution to the Ricci flow starting from the compact approximating manifold $(M_i,g_{i,0})$ obtained from Proposition~\ref{prop:CW-appr}, then 
 $$\lim_{i\to+\infty} T_i=+\infty.$$
\end{prop}
\begin{proof}
Let $\hat T>0$, we claim that $T_i>\hat T$ for all $i$ sufficiently large. Suppose on the contrary that $T_{i_k}\leq \hat T$ for some $i_k\to+\infty$. For convenience, we still denote it as $i$. Lemma~\ref{eqn:evo-W} implies 
\begin{equation}\label{eqn:upp-bdd-W}
\begin{split}
  W_{i,R,m}(t)&\leq \left(C_1m^2 R^{-2+\frac{2(n-2)}{m}} t+\left(W_{i,R,m}(0)\right)^\frac2m\right)^\frac{m}2\\
  &\leq C_2\left(C_1^\frac{m}2m^2 R^{-m+n-2} t^\frac{m}2+W_{i,R,m}(0)\right)\\
    &\leq C_2\left(C_1^\frac{m}2m^2 R^{-m+n-2} \hat T^\frac{m}2+W_{i,R,m}(0)\right)=:\mathbf{W}_{i,R,0}(\hat T)
  \end{split}
\end{equation}
for $t\in [0,T_i)$, where $C_2>1$ is a constant depending only on $m$.

By integrating the evolution inequality in Lemma~\ref{eqn:evo-V} and using  \eqref{eqn:upp-bdd-W}, we conclude that for $t\in [0,T_i)\subseteq [0,\hat T)$,
\begin{equation}\label{eqn:esti-V-W-etc}
\begin{split}
\mathrm{Vol}_{g_i(t)}(M_i)&\geq V_{i,R,m}(t)\\
&\geq V_{i,R,m}(0)-\mathbf{W}_{i,R,0} (\hat T)\cdot \hat T -\frac{10^33^n  m^2\hat T}{R^2} \cdot\mathrm{Vol}_{g_{i,0}}\left(p_i,R\right).
\end{split}
\end{equation}

Since the right hand side is independent of $t\in [0,T_i)$, we 
 use Lemma~\ref{lma:lifespan} and let $t\to T_i$ to conclude 
\begin{equation}\label{eqn:esti-V-W-etc-1}
\begin{split}
V_{i,R,m}(0)\leq \mathbf{W}_{i,R,0} (\hat T)\cdot \hat T +\frac{10^33^n  m^2\hat T}{R^2} \cdot\mathrm{Vol}_{g_{i,0}}\left(p_i,R\right)
\end{split}
\end{equation}

From Lemma~\ref{lma:vol-init} and Lemma~\ref{lma:init-W}, for $R$ sufficiently large we have 
\begin{equation}\label{eqn:initi-V--}
\lim_{i\to+\infty} V_{i,R,m}(0)\geq 2^{-n}\mathrm{Vol}_{g_0}(p_0,R)
\end{equation}
and 
\begin{equation}\label{eqn:initi-W--}
\begin{split}
\lim_{i\to +\infty} \mathbf{W}_{i,R,0} (\hat T)&\leq C_2 4^{n-2}R^{-2} k(p_0,4R)\cdot \mathrm{Vol}_{g_0}(p_0,R)\\
&\quad +C_2C_1^{m/2}m^2 R^{-m+n-2}\hat T^{m/2}.
\end{split}
\end{equation}

Together with the local smooth convergence of $g_{i,0}$ to $g_0$ from Proposition~\ref{prop:CW-appr},  \eqref{eqn:esti-V-W-etc-1}, \eqref{eqn:initi-V--} and \eqref{eqn:initi-W--} imply
\begin{equation}\label{eqn:V-W--interact}
\begin{split}
&\quad 2^{-n}\cdot\mathrm{Vol}_{g_0}(p_0,R)\\
&\leq  C_24^{n-2}\hat TR^{-2} k(p_0,4R)\cdot \mathrm{Vol}_{g_0}(p_0,R)\\
&\quad +C_2C_1^{m/2}m^2 R^{-m+n-2}\hat T^{1+\frac{m}2}+\frac{10^33^n  m^2\hat T}{R^2} \cdot\mathrm{Vol}_{g_{0}}\left(p_0,R\right).
\end{split}
\end{equation}

\medskip

We now fix $m>n$ so that the exponent of $R$ appeared in the second term on the right hand side is negative, i.e. $-m+n-2\leq -2$. Since $\mathrm{Vol}_{g_0}(p_0,R)\geq  \mathrm{Vol}_{g_0}(p_0,1):=v_0>0$ and  $R^{-2} k(p_0,4R)=o(1)$ as $R\to+\infty$, we have contradiction by passing $R\to+\infty$. This completes the proof.
\end{proof}

\subsection{Uniform local regularity \texorpdfstring{$g_i(t)$}{gi(t)}}

From Proposition~\ref{prop:almost-immort}, we showed that $g_i(t)$ is almost immortal as $i\to+\infty$. We next show that $g_i(t)$ converges locally uniformly smoothly as $i\to+\infty$. In particular, we construct an immortal solution to the Ricci flow $g(t)$ starting from $(M,g_0)$. For our purpose, we will construct Ricci flow with Li-Yau-Hamilton trace Harnack inequality. We say that a complete solution $(M,g(t))$, $t\in (0,T)$, to the Ricci flow 
satisfies Li-Yau-Hamilton trace Harnack inequality if 
for all $(p,t)\in M\times (0,T)$ and $V\in T_pM$, 
\[
\frac{\partial}{\partial t}\mathrm{scal}+\frac{1}{t}\mathrm{scal}+2\langle \nabla \mathrm{scal}, V\rangle+2\Ric(V,V)\ge 0.
\] 

In \cite{Brendle2009}, Brendle showed that a complete Ricci flow with $\mathrm{sec}^{\mathbb{C}}(g(t))\geq 0$ and uniformly bounded curvature on compact time intervals must satisfy the  Li-Yau-Hamilton matrix Harnack inequality, in particular the trace Harnack inequality,  see also \cite{Hamilton1993, Cao1992} for earlier work. We now use it to construct a global solution under some average curvature decay.

\begin{thm}\label{thm:immort}
Suppose $(M^n,g_0)$ is a complete non-compact manifold with $\mathrm{sec}^\mathbb{C}\geq 0$  and 
\begin{equation}
 \lim_{r\to+\infty}r^{-2}k(x_0,r)=0,\;\;\text{as}\;r\to+\infty
\end{equation}
for some $x_0\in M$, then it admits a Ricci flow $g(t)$ on $\tilde M\times [0,+\infty)$ such that $g(0)=\tilde g_0$, $\mathrm{sec}^\mathbb{C}(g(t))\geq 0$ and satisfies the Li-Yau-Hamilton Harnack inequality, here $(\tilde M, \tilde g_0)$ is the universal cover of $(M,g_0)$. 
\end{thm}
\begin{proof}
By Lemma~\ref{lma:univer-k}, $(\tilde M,\tilde g_0)$ also satisfies the same average curvature decay assumption for some $\tilde x_0\in \tilde M$. Furthermore $\tilde M^n$ is diffeomorphic to $\mathbb{R}^n$ by Proposition~\ref{prop:soul-pt}. In particular, its soul is a point $p_0$. By Lemma~\ref{lma:indep-centr}, we also have 
\begin{equation}
k(p_0,r)=o(r^2)
\end{equation}
as $r\to+\infty$. 
Therefore without loss of generality, by passing to the universal covering if necessary, we might assume $M=\tilde M$, $x_0=\tilde x_0=p_0$ is a  soul point and $M$ is diffeomorphic to $\mathbb{R}^n$.

We now claim that the Ricci flow constructed in \cite{CabezasWilking2015} is in fact immortal. This is constructed by taking limit of $g_i(t)$ where $(M_i,g_i(t)),t\in [0,T_i)$ is the sequence of Ricci flows in Proposition~\ref{prop:almost-immort}.  In \cite{CabezasWilking2015}, it has already showed that $g_i(t)$ is locally uniformly $C^\infty$ regular for $t\in [0,T]$ for some uniform small $T>0$. It remains to give estimate for $t>T$.

By Brendle's generalization of (trace) Li-Yau-Hamilton's Harnack inequality \cite{Brendle2009}, for any $x_1,x_2\in M_i$ and $0<t_1<t_2$ we have 
\begin{equation}\label{eqn:harn-}
\begin{split}
\mathrm{scal}(g_i(x_1,t_1))&\leq \frac{t_2}{t_1}\mathrm{scal}(g_i(x_2,t_2))\cdot \exp\left(\frac{d_{g_{i,0}}(x_1,x_2)^2}{4(t_2-t_1)} \right).
\end{split}
\end{equation}
where we have also used $g_i(t)\leq g_{i,0}$ and $2\Ric (g_i(t))\leq \mathrm{scal}(g_i(t))$ for $t\in [0,T_i)$.

\begin{claim}\label{claim:harnk-int}
For all $t>0$, $s\in (t,T_i)$, sufficiently large $R>0$ and $i,m\in \mathbb{N}$, we have 
\begin{equation}
\begin{split}
 \mathrm{scal}(g_i(p_i,t))\cdot V_{i,R,m}(s)
&\leq \exp\left(\frac{4R^2}{s-t} \right)\cdot \frac{s}{t}\cdot  W_{i,R,m}(s).
\end{split}
\end{equation}
\end{claim}
\begin{proof}[Proof of Claim~\ref{claim:harnk-int}]

By \eqref{eqn:harn-}, if we fix $x_1=p_i$ and $z=x_2\in B_{g_{i,0}}(p_i,4R)$, then 
\begin{equation}
\begin{split}
\mathrm{scal}(g_i(p_i,t))\leq \frac{s}{t}\cdot \mathrm{scal}(g_i(z,s))\cdot \exp\left(\frac{4R^2}{s-t} \right)
\end{split}
\end{equation}
for all $0<t<s<T_i$.

We multiply both sides by $\Phi_{i,R,m}(z,s)$ and integrate it over $B_{g_{i,0}}(p_i,4R)$ at time $s$ to conclude 
\begin{equation}\label{eqn:V-Harn-W}
\begin{split}
&\quad \exp\left(-\frac{4R^2}{s-t} \right)\cdot \frac{t}{s}\mathrm{scal}(g_i(p_i,t)) \int_{B_{g_{i,0}}(p_i,4R)} \Phi_{i,R,m}(z,s)\,d\mathrm{vol}_{g_i(z,s)}\\
&\leq \int_{B_{g_{i,0}}(p_i,4R)}\mathrm{scal}(g_i(z,s)) \Phi_{i,R,m}(z,s)\,d\mathrm{vol}_{g_i(z,s)}\\
&\leq W_{i,R,m}(s).
\end{split}
\end{equation}

We rewrite the integral on the left hand side by:
\begin{equation}\label{eqn:int-cut}
\begin{split}
&\quad \int_{B_{g_{i,0}}(p_i,4R)} \Phi_{i,R,m}(z,s)\,d\mathrm{vol}_{g_i(z,s)}\\
&:=V_{i,R,m}(s)-\int_{M_i\setminus B_{g_{i,0}}(p_i,4R)} \Phi_{i,R,m}(z,s)\,d\mathrm{vol}_{g_i(z,s)}.
\end{split}
\end{equation}
If $z\in M_i$ is such that $d_{g_{i,0}}(p_i,z)\geq 4R$, then 
\begin{equation}
u_i(z,s)\geq b_i(z)\geq \frac12 d_{g_{i,0}}(p_i,z)\geq 2R
\end{equation}
and hence $\Phi_{i,R,m}(z,s)=0$ so that the last term in \eqref{eqn:int-cut} vanishes. This proves the claim by \eqref{eqn:V-Harn-W}.
\end{proof}

Following Claim~\ref{claim:harnk-int},  we now estimate $W_{i,R,m}(s)$ from above and $V_{i,R,m}(s)$ from below.  By \eqref{eqn:esti-V-W-etc} with $\hat T$ replaced by $s\in (0,T_i)$,  we have 
\begin{equation}\label{eqn:VVVV}
\begin{split}
 V_{i,R,m}(s)
&\geq V_{i,R,m}(0)-\mathbf{W}_{i,R,0} (s)\cdot s  -\frac{10^33^n  m^2s}{R^2} \cdot\mathrm{Vol}_{g_{i,0}}\left(p_i,R\right)
\end{split}
\end{equation}
where 
\begin{equation}
\mathbf{W}_{i,R,0} (s)=C_2C_1^\frac{m}2m^2 R^{-m+n-2} s^\frac{m}2+C_2W_{i,R,m}(0).
\end{equation}

By \eqref{eqn:upp-bdd-W}, \eqref{eqn:initi-V--} and \eqref{eqn:initi-W--}, for $R$ sufficiently large and $m\in \mathbb{N}$ we choose $N(R,m)$ such that for $i>N$, we have
\begin{equation} \label{eqn:inii}
\begin{split}
V_{i,R,m}(0)&> 2^{-n-1}\mathrm{Vol}_{g_0}(p_0,R)\\[2mm]
 W_{i,R,m}(s)&\leq \mathbf{W}_{i,R,0} (s)< C_24^{n-1}R^{-2} k(p_0,4R)\cdot \mathrm{Vol}_{g_0}(p_0,R)\\
&\quad +2C_2C_1^{m/2}m^2 R^{-m+n-2}s^{m/2}.
\end{split}
\end{equation}

By putting \eqref{eqn:inii} into \eqref{eqn:VVVV}, we conclude for $R$ sufficiently large, $i>N$ and $s\in (0,T_i)$ that 
\begin{equation}\label{eqn:V_i(s)}
\begin{split}
V_{i,R,m}(s)&\geq \mathrm{Vol}_{g_0}(p_0,R)\cdot \left(\frac1{2^{n+2}}-\frac{4^{n-1}C_2s k(p_0,4R)}{R^2} -\frac{10^33^n m^2 s}{R^2}  \right)\\
&\quad  +\left(\frac1{2^{n+2}}\mathrm{Vol}_{g_0}(p_0,R)-2C_2C_1^{m/2} m^2 R^{-m-2+n}s^{1+\frac{m}2}\right).
\end{split}
\end{equation}

For any large $\hat T>0$, since $T_i\to+\infty$ from Proposition~\ref{prop:almost-immort}, we may suppose $T_i>\hat T$ for all large $i$, as $R\to+\infty$ if we restrict $s\in (0,\hat T]$ and choose $m=4n+4$, then \eqref{eqn:V_i(s)} and $k(p_0,R)=o(R^2)$ shows that for $R$ sufficiently  large (depending on $\hat T$ but uniform in $i$),  $i>N$ and $s\in (0,\hat T]$,
\begin{equation}\label{eqn:V-asy}
V_{i,R,m}(s)\geq  2^{-n-3}\mathrm{Vol}_{g_0}(p_0,R).
\end{equation}
Combines \eqref{eqn:V-asy} with \eqref{eqn:inii} and Claim~\ref{claim:harnk-int}, we conclude that for all such large $R>0$ there exists $C_R>0$ such that for all $i\to+\infty$ and $t\in (0,\hat T]$,
\begin{equation}
\mathrm{scal}(g_i(p_i,t))\leq C_Rt^{-1}.
\end{equation}
By $\mathrm{sec}^\mathbb{C}(g_i(t))\geq 0$, Harnack inequality \eqref{eqn:harn-}, Shi's estimate \cite{Shi1989} and Proposition~\ref{prop:almost-immort}, we conclude that the Ricci flow constructed in \cite{CabezasWilking2015} is bounded in $C^k_{loc}$ locally uniformly on $M\times [0,+\infty)$  
and hence by a diagonal sequence argument, it converges locally uniformly smoothly to an immortal solution $g(t)$ on $M\times [0,+\infty)$. The completeness of $g(t)$ for general $t>0$ follows from \cite[Lemma 4.9]{CabezasWilking2015}. Since $g(t)$ is local limit of compact Ricci flow with $\mathrm{sec}^\mathbb{C}\geq 0$, Harnack inequality follows directly from smooth convergence.
\end{proof}

\begin{rem}
    We note here that $g(t)$ constructed in Theorem~\ref{thm:immort} does not necessarily have bounded curvature for $t>0$. Nevertheless, it satisfies Li-Yau-Hamilton Harnack inequality by local smooth convergence of compact Ricci flow $g_i(t)$ with $\mathrm{sec}^\mathbb{C}(g_i(t))\geq 0$. 
\end{rem}

\medskip
Finally, we note that the function $u_i$ in Lemma~\ref{lma:heat-flowed} sub-converges to a function $u$ on $M\times [0,+\infty)$ with the same properties. 
\begin{prop}\label{prop:v-ch-fnc}
The function $u_i$ in Lemma~\ref{lma:heat-flowed} sub-converges to $v$ locally uniformly in $M\times [0,+\infty)$ such that 
\begin{enumerate}
    \item [(i)] $v$ is $1$-Lipschitz;
    \item[(ii)] $0\leq \partial_t v\leq \Delta_{g(t)}v$ in distribution sense;
    \item[(iii)] $v(x,t)\geq v(x,0)=\b(x)$ for $x\in M$.
\end{enumerate}
\end{prop}
\begin{proof}
Since $u_i$ is $1$-Lipschitz with respect to $g_i(t)$ for each $t\in [0,+\infty)$ and $g_i(t)$ converges to $g(t)$ locally uniformly smoothly, it suffices to show that its time derivatives is locally uniformly bounded. (i) and (ii) will then follow from Arzel\`a-Ascoli theorem and Lemma \ref{lma:heat-flowed}.
We fix $\Omega\Subset M,S>0$ and choose $A,\rho>0$ such that for all $i\to+\infty$, 
\begin{equation}
    0\leq \Ric(g_i(t))\leq A g_i(t)
\end{equation}
on $B_{g_i(t)}(x,\rho)$, for all $x\in \Omega$ and $t\in [0,S]$. Let $\gamma:[0,L]\to M_i$ be the normal geodesic from $x=\gamma(0)$ to $N_i$, then 
\begin{equation}
    \partial_t u_i(x,t)=\int_0^{L} \Ric_{g_i(t)}(\gamma',\gamma')ds.
\end{equation}

As in the proof of \cite[Lemma 8.3]{Perelman} and \cite[Lemma 4.8]{CabezasWilking2015}, we use the fact that $N_i$ is totally geodesic so that second variation applies and shows that $\partial_t u_i(x,t)\leq (n-1)\rho^{-1}+A\rho$. $v(x,0)=\b(x)$ follows from Proposition \ref{prop:CW-appr} (v), Lemma \ref{lma:heat-flowed} and the fact that $u_i(x,0)=b_i(x)$.
\end{proof}

\section{Volume non-collapsing by local smoothing}

In this section, we consider the scaling invariant case of Theorem~\ref{thm:immort}. That said we consider the case  when $k(x_0,r)$ is uniformly bounded for all $r>0$. This follows from modifying the weight used in section~\ref{sec:long-time}.

\begin{thm}\label{thm:soul}
Suppose $g$ is a complete metric on $\mathbb{R}^n$ 
such that 
\begin{enumerate}
\item[(ii)] $\mathrm{sec}^\mathbb{C}\geq 0$;
\item[(iii)] there is $p_0\in M,\Lambda>0$ such that $\sup_{r>0}k(p_0,r)\leq \Lambda$,
\end{enumerate}
then 
$$\mathrm{AVR}(g):=\limsup_{r\to+\infty}\frac{\mathrm{Vol}_g(x_0,r)}{\omega_nr^n}>0.$$
\end{thm}
\begin{proof}

We let $g(t)$ be the immortal Ricci flow constructed in Theorem~\ref{thm:immort} and $v(x,t)$ be the heat solution starting from $\b(x)$ (in $L^\infty$ sense) from Proposition~\ref{prop:v-ch-fnc}. We will use $C_i$ to denote any dimensional constants. We first show that the Ricci flow $g(t)$ indeed satisfies the scaling invariant curvature estimates $\sup_M|\Rm(g(t))|=O(t^{-1})$ as $t\to+\infty$ and $\mathrm{sec}^\mathbb{C}(g(t))\geq 0$.

We modify the choice of cut-off function $\Phi_{\infty,R,m}$ used in \eqref{eqn:defining-W-V}. For $R>0$, we define 
\begin{equation}
\Phi_R(x,t):=\phi\left(\frac{v(x,t)}{R}\right)=
\exp\left(-\frac{v(x,t)}{R} \right)
\end{equation}
and consider 
\begin{equation}\label{eqn:defining-W-V-limit}
\left\{
\begin{array}{ll}
\displaystyle V_{R}(t):=\int_{M} \Phi_{R}\,d\mathrm{vol}_{g(t)}\\[3mm]
\displaystyle W_{R}(t):=\int_{M} \Phi_{R}\cdot\mathrm{scal}(g(t))\,d\mathrm{vol}_{g(t)}
\end{array}
\right.
\end{equation}
for $t\geq 0$. We first see that both $V_R$ and $W_R$ are finite for any $t\geq 0$.

\begin{claim}\label{claim:int}
For every $t\in [0,+\infty)$, the function $\Phi_R$ and $\Phi_R \cdot\mathrm{scal}(g)$ is integrable with respect to $d\mathrm{vol}_{g(t)}$. Moreover, there exists $L_1(n)>0$ such that $$V_R(t)\leq L_1  \mathrm{Vol}_{g(t)}(p_0,R)\quad\text{and}\quad W_R(t)\leq L_1R^{n-2}$$ for all $t\in [0,+\infty)$. Furthermore 
\end{claim}
\begin{proof}[Proof of Claim~\ref{claim:int}]

Recall that for all $(x,t)\in M\times [0,+\infty)$, we have 
\begin{equation}
v(x,t)\geq \b(x)\geq \frac12 d_{g_0}(x,p_0)\geq \frac12 d_{g(t)}(x,p_0)
\end{equation}
and hence co-area formula and volume comparison imply
\begin{equation}\label{eqn:Vol_R(t)}
\begin{split}
V_R(t)&\leq \int_{M\setminus B_{g(t)}(p_0,R)} \exp\left(-\frac{d_{g(t)}(p_0,\cdot)}{2R} \right)\,d\mathrm{vol}_{g(t)}+ \mathrm{Vol}_{g(t)}(p_0,R)\\
&\leq \int^{+\infty}_R \exp\left(-\frac{r}{2R} \right) |\partial B_{g(t)}(p_0,r)|\,dr + \mathrm{Vol}_{g(t)}(p_0,R)\\
&\leq C_0  \mathrm{Vol}_{g(t)}(p_0,R) .
\end{split}
\end{equation}

For $W_R$, we use co-area formula, volume comparison and Stokes' Theorem to obtain
\begin{equation}\label{eqn:W_R-finite}
\begin{split}
W_R(t)&\leq \int_M \mathrm{scal}(g(t))\cdot \exp\left(-\frac{d_{g(t)}(p_0,\cdot)}{2R} \right)\,d\mathrm{vol}_{g(t)}\\
&=\int^{+\infty}_0 \exp\left(-\frac{r}{2R} \right)\cdot \left(\int_{\partial B_{g(t)}(p_0,r)} \mathrm{scal}(g(t)) dA_{g(t)}\right)\,dr\\
&= \exp\left(-\frac{r}{2R} \right)\cdot\left(\int_{ B_{g(t)}(p_0,r)} \mathrm{scal}(g(t)) d\mathrm{vol}_{g(t)}\right)\Big|_{r=0}^{r=+\infty}\\
&\quad +\frac1{2R}\int^{+\infty}_0 \exp\left(-\frac{r}{2R} \right)\cdot\left(\int_{ B_{g(t)}(p_0,r)} \mathrm{scal}(g(t)) d\mathrm{vol}_{g(t)}\right)dr.
\end{split}
\end{equation}

By the rescaled version of  Petrunin's estimate \cite{Petrunin2009}, see \cite[Lemma 7.1]{CabezasWilking2015}, 
\begin{equation}
\int_{ B_{g(t)}(p_0,r)} \mathrm{scal}(g(t)) d\mathrm{vol}_{g(t)}\leq C_1r^{n-2}
\end{equation}
Substituting it into \eqref{eqn:W_R-finite} yields $
W_R(t)\leq C_2 R^{n-2}$ for all $t\in [0,+\infty)$. 

\end{proof}

Next we have the analogy to Lemma~\ref{lma:vol-init} and Lemma~\ref{lma:init-W}.
\begin{claim}\label{claim:init-limit}
There is $L_2(n)>0$ such that  
\begin{equation}
\left\{
\begin{array}{ll}
V_R(0)\geq L_2^{-1} \mathrm{Vol}_{g_0}(p_0,R);\\[3mm]
W_R(0)\leq L_2 \Lambda R^{-2} \mathrm{Vol}_{g_0}(p_0,R).
\end{array}
\right.
\end{equation}
\end{claim}

\begin{proof}[Proof of Claim~\ref{claim:init-limit}]
The estimate of $V_R(0)$ follows the same argument in  Lemma~\ref{lma:vol-init}.

For $W_R(0)$,  we follow the derivation in \eqref{eqn:W_R-finite} except we use $k(p_0,r)$ instead:
\begin{equation}
\begin{split}
W_R(0)&\leq \int_{B_{g_0}(p_0,R)}\mathrm{scal}(g_0)\,d\mathrm{vol}_{g_0}+\int_{M\setminus B_{g_0}(p_0,R)}\Phi_R\cdot \mathrm{scal}(g_0)\,d\mathrm{vol}_{g_0}\\
&\leq \frac{\Lambda}{R^2}\cdot \mathrm{Vol}_{g_0}(p_0,R)+\int^{+\infty}_R\exp\left(-\frac{r}{2R} \right)\cdot \left(\int_{\partial B_{g(t)}(p_0,r)} \mathrm{scal}(g(t)) dA_{g(t)}\right)\,dr\\
&\leq \frac{\Lambda}{R^2}\cdot \mathrm{Vol}_{g_0}(p_0,R)+\frac1{2R}\int^{+\infty}_R \exp\left(-\frac{r}{2R} \right)\cdot\left(\int_{ B_{g_0}(p_0,r)} \mathrm{scal}(g_0) d\mathrm{vol}_{g_0}\right)dr\\
&\leq  \frac{\Lambda}{R^2}\cdot \mathrm{Vol}_{g_0}(p_0,R)+\frac{\Lambda}{2R}\int^{+\infty}_R \exp\left(-\frac{r}{2R} \right)\cdot \frac{1}{r^2}\cdot \mathrm{Vol}_{g_0}(p_0,r)dr\\
&\leq   \frac{\Lambda}{R^2}\cdot \mathrm{Vol}_{g_0}(p_0,R)+\frac{\Lambda \mathrm{Vol}_{g_0}(p_0,R)}{2 R^{n+1}}\int^{+\infty}_R \exp\left(-\frac{r}{2R} \right)\cdot r^{n-2}\cdot dr\\
&\leq C_3\Lambda R^{-2}\mathrm{Vol}_{g_0}(p_0,R),
\end{split}
\end{equation}
where we have used the volume comparison in the second last inequality.
\end{proof}

The major advantage of using exponential decay $\Phi_R$ instead of cut-off function $\Phi_{i,R,m}$ is to improve the evolution inequality in Lemma~\ref{eqn:evo-V} and Lemma~\ref{eqn:evo-W}, using the fact that $\phi''=\phi$ and $\phi'\leq 0$. With the new choice of weight, we have the following  inequalities directly modified from Lemma~\ref{eqn:evo-V} and Lemma~\ref{eqn:evo-W}: 
\begin{equation}\label{eqn:mod-evo-WV}
\left\{
\begin{array}{ll}
\displaystyle \partial_t V_R\geq -W_R-\frac1{R^2}V_R;\\[3mm]
\displaystyle \partial_t W_R\leq  \frac1{R^2}W_R.
\end{array}
\right.
\end{equation}
Here the derivatives is understood in sense of Dini derivatives. The Stoke Theorem involved can be justified using the exponential decay as in the proof of Claim~\ref{claim:int}.

From \eqref{eqn:mod-evo-WV} and Claim~\ref{claim:init-limit},  
\begin{equation}\label{eqn:improved-W}
\begin{split}
W_R(t)&\leq W_R(0)\cdot \exp\left(\frac{t}{R^2} \right)\\
&\leq  3L_2\Lambda R^{-2}  \cdot\mathrm{Vol}_{g_0}(p_0,R).
\end{split}
\end{equation}
for $t\in [0,R^2]$.

This improves the estimate \eqref{eqn:upp-bdd-W} where the weight is cut-off. Consequently, we use \eqref{eqn:improved-W}, \eqref{eqn:mod-evo-WV} and Claim~\ref{claim:init-limit} to deduce
\begin{equation}\label{eqn:low-V_R(t)}
\begin{split}
3 V_R(t)&\geq V_R(0)-10 L_2 \Lambda t R^{-2} \cdot\mathrm{Vol}_{g_0}(p_0,R)\\
&\geq \left(L_2^{-1}-10 L_2 \Lambda t R^{-2} \right)\cdot \mathrm{Vol}_{g_0}(p_0,R)\\
&\geq \frac1{2L_2} \mathrm{Vol}_{g_0}(p_0,R)\\
\end{split}
\end{equation}
for $t\in [0,\delta R^2]$ where $\delta:=\Lambda^{-1}(10L_2)^{-2}<1$. We will fix this choice of $\delta$ from now.

\medskip

Finally we need the analogy to Claim~\ref{claim:harnk-int}. Similar to the derivation of \eqref{eqn:V-Harn-W}, we fix a large $\tilde \Lambda$ and integrate the Harnack estimate over $B_{g_0}(p_0,2\tilde \Lambda R)$ so that for all $0<t<s$, 
\begin{equation}\label{eqn:V-Harn-W-limit}
\begin{split}
 \exp\left(-\frac{\tilde\Lambda^2 R^2}{s-t} \right)\cdot \frac{t}{s}\mathrm{scal}(g(p_0,t)) \int_{B_{g_{0}}(p_0,2\tilde\Lambda R)} \Phi_{R}(z,s)\,d\mathrm{vol}_{g(z,s)}
&\leq W_{R}(s)
\end{split}
\end{equation}
where $$v(x,s)\geq \b(x)\geq \frac12 d_{g_0}(x,p_0),$$
 the same argument as \eqref{eqn:Vol_R(t)} implies
\begin{equation}\label{eqn:int-cut-2}
\begin{split}
\int_{M\setminus B_{g_{0}}(p_0,2\tilde\Lambda R)} \Phi_{R}(z,s)\,d\mathrm{vol}_{g(z,s)}
&\leq e^{-\tilde\Lambda /2}\cdot\int_M \exp\left(-\frac{d_{g_0}(\cdot ,p_0)}{4R} \right)d\mathrm{vol}_{g_0}\\
&\leq C_4 e^{-\tilde\Lambda/2} \cdot\mathrm{Vol}_{g_0}(p_0,R)\\
&\leq  6C_4 L_2 e^{-\tilde\Lambda/2} V_R(s)
\end{split}
\end{equation}
if $s\in [0,\delta R^2]$, where we have used \eqref{eqn:low-V_R(t)}.

By substituting \eqref{eqn:int-cut-2} into \eqref{eqn:V-Harn-W-limit}, we use  \eqref{eqn:improved-W} to deduce
\begin{equation}
\begin{split}
&\quad  \exp\left(-\frac{\tilde\Lambda^2 R^2}{s-t} \right)\cdot \frac{t}{s}\mathrm{scal}(g(p_0,t)) \cdot\left(1-6C_4L_2 e^{-\tilde \Lambda/2} \right) V_R(s)\\
&\leq 20L_2^2\Lambda R^{-2} \cdot V_R(s).
\end{split}
\end{equation}
for $0<t<s\leq \delta R^2$. 

We first choose $\tilde \Lambda(n)$ large enough so that $6C_4L_2 e^{-\tilde \Lambda/2}\leq 1/2$. For $R\to+\infty$, we choose $s=s_R=\delta R^2$ and $t=t_R=s_R/2$ to conclude that 
\begin{equation}
\mathrm{scal}(p_0,\delta R^2/2)\leq \a(n,\Lambda)R^{-2}
\end{equation}
for all $R\to+\infty$. By applying Li-Yau-Hamilton's Harnack inequality \cite{Brendle2009} again, we deduce that 
for $(x,t)\in B_{g_0}(p_0,R)\times (0,\delta R^2/4)$
\[
\mathrm{scal}(x,t)\leq \frac{\delta\a}{2t}\exp\left(\frac{d_{g_0}(x,p_0)^2}{\delta R^2}\right)\leq \frac{\delta\a}{2t}\exp\left(\delta^{-1}\right).
\]
Hence
$$\sup_M|\Rm(g(t))|\leq \a_1 t^{-1}$$
for all $t\in (0,+\infty)$. 

We now use the strong maximum principle to reduce the discussion to $\mathrm{sec}^\mathbb{C}>0$, following \cite[Theorem 8.1]{DeruelleLeeSchulzeSimomTopping2026} and \cite[Corollary 7.6]{CabezasWilking2015}. We fix $s_0$ small and consider the de-Rham
decomposition of $(M,g(s_0))$:
\begin{equation}
(M,g(s_0))=\mathbb{R}^k\times \prod_{i=1}^N (N_i,h_i)
\end{equation}
where each $(N_i,h_i)$ is irreducible. By the existence and uniqueness \cites{ChenZhu2006,Shi1989} of Ricci flow with bounded curvature, the Ricci flow $(M,g(t))$ splits as $\mathbb{R}^k\times \prod_{i=1}^N (N_i,h_i(t))$ for all $t\geq s_0$ and each $h_i(t),t\in [s_0,+\infty)$ satisfies 
\begin{equation}
0\leq \mathrm{sec}^\mathbb{C}(h_i(t))\leq \a t^{-1}
\end{equation}
for all $i$ and $t\in [s_0,+\infty)$. Moreover by the strong maximum principle, we might assume $\mathrm{scal}(h_i(t))>0$ on $(s_0,+\infty)$ since otherwise $\Rm(h_i(t))\equiv 0$ which is impossible. We might assume $h_i(t)$ is irreducible for $t\approx s_0$. We claim that $\mathrm{AVR}(h_i(t))>0$ for $t>s_0$.

 Arguing as in \cite[Theorem 8.1]{DeruelleLeeSchulzeSimomTopping2026}, since $(N_i,h_i(t)),t\in (s_0,+\infty)$ has $\mathrm{scal}(h_i(t))>0$, it cannot be Ricci flat. Moreover it must be non-symmetric, non-compact and non-Einstein, since otherwise the positivity of scalar curvature implies a finite time singularity which is impossible.  Now Berger’s holonomy classification theorem implies that the holonomy group must be either $\mathrm{SO}(n)$ or $\mathrm{U}(n/2)$. In the first case, it follows from \cite[Proposition 9]{BrendleSchoen2008} that $\mathrm{sec}^\mathbb{C}(h_i(t))>0$ for $t>s_0$ since it is non-flat. In the latter case, we might argue similarly that at least we have  $\mathrm{sec}(h_i(t))>0$ for $t>s_0$, for example see the proof of \cite[Theorem 4.4]{LeeTopping2025} or \cite[Corollary 7.6]{CabezasWilking2015}.  

Now we are in position to carry out the argument in the proof of \cite[Theorem 4.4]{LeeTopping2025}. We include it for readers' convenience. By a result of Gromoll-Meyer \cite[Theorem B.65]{ChowKnop},  there is $c_0>0$ such that $\mathrm{inj}(h_i(t))\geq c_0\sqrt{t}$ for all $t>s_0$ and $i$. Since $s_0$ is arbitrary, this shows that $\mathrm{inj}(g(t))\geq c_0'\sqrt{t}$ for all $t>0$. Therefore a Theorem of G\"unther implies 
\begin{equation}
\mathrm{Vol}_{g(t)}(x_0,c_0\sqrt{t}) \geq c_1 t^{n/2}
\end{equation}
for some $c_1>0$ and any $x_0\in M$. Thanks to the curvature decay, \cite[Corollary 3.3]{SimonToppingJDG} implies that for some $c_2>0$ and for all $t>0$,
\begin{equation}
\mathrm{Vol}_{g_0}(x_0,c_2\sqrt{t})\geq \mathrm{Vol}_{g(t)}(x_0,c_0\sqrt{t})\geq c_1t^{n/2}.
\end{equation}
Since $t$ is arbitrarily large, this shows that $\mathrm{AVR}(g_0)>0$. Alternatively one can also prove using the method of \cite[Theorem 2.2]{ChauTamTAMS}. This completes the proof.
\end{proof}

We are now ready to give a proof to our main results, i.e. Theorem~\ref{thm:main} and Corollary~\ref{cor:gap}.
\begin{proof}[Proof of Theorem~\ref{thm:main}]

By Lemma~\ref{lma:univer-k} and Proposition~\ref{prop:soul-pt}, the universal cover $(\tilde M,\tilde g)$ is diffeomorphic to $\mathbb{R}^n$ and satisfies 
\begin{equation}
\sup_{r>0} k_{\tilde g}(\tilde x_0,r)<+\infty
\end{equation}
for some $\tilde x_0\in \tilde M$. Hence Theorem~\ref{thm:soul} implies that  $\mathrm{AVR}(\tilde g)>0$. 
\end{proof}

With the non-collapsing from Theorem~\ref{thm:main}, the gap Theorem is immediate. 
\begin{proof}[Proof of Corollary~\ref{cor:gap}]

By Lemma~\ref{lma:univer-k} and Proposition~\ref{prop:soul-pt}, the universal cover $(\tilde M,\tilde g)$ satisfies $\mathrm{AVR}(\tilde g)>0$ and 
\begin{equation}
k_{\tilde g}(\tilde x_0,r)=o(1)
\end{equation}
as $r\to+\infty$, for some $\tilde x_0\in \tilde M$.   By Theorem~\ref{thm:main}, we have $\mathrm{AVR}(\tilde g)>0$. The flatness follows from \cite[Theorem 1.3]{ChanLee2025}.
\end{proof}

\end{document}